\input amstex
\documentstyle{amsppt}
\output={\plainoutput}
\pageno=1
\footline={\ifnum\pageno>1 \myfl\fi}
\define\myfl{\hss\tenrm\folio\hss}
\headline={\hfil}
\NoBlackBoxes
\magnification=\magstep 1
\pagewidth{15 true cm}
\pageheight{22 true cm}

\topmatter
\title Burghelea-Friedlander-Kappeler's gluing formula and the
adiabatic decomposition of the zeta-determinant of a Dirac Laplacian.
\endtitle
\author Yoonweon Lee
\endauthor
\affil  Department of Mathematics \\
        Inha University \\
        Inchon, 402-751, Korea
\endaffil
\endsubjclass
\thanks
This work was supported partially by 1999-Inha research fund.
\endthanks
\abstract
In this paper we first establish the relation between the zeta-determinant of a Dirac Laplacian with the
Dirichlet boundary condition and the APS boundary condition on a cylinder. Using this result and the gluing formula of the zeta-determinant given by Burghelea, Friedlander and Kappeler with some assumptions,
we prove the adiabatic decomposition theorem of the zeta-determinant of a Dirac Laplacian. This result was originally proved by J. Park and 
K. Wojciechowski in [11] but our method is completely different from the one they presented.
\endabstract
\endtopmatter

\document
\baselineskip 0.6 true cm plus 3 pt minus 3 pt

\S 1 {\bf Introduction}
\TagsOnRight
\vskip 0.5 true cm

Let $M$ be a compact oriented $m$-dimensional Riemannian manifold and $E\rightarrow M$ be a Clifford module bundle.
Suppose that ${\frak D}_{M}$ is a Dirac operator acting on smooth sections of $E$ .
Then ${\frak D}_{M}^{2}$ is a non-negative self-adjoint elliptic differential operator, which is called a Dirac Laplacian.
From the standard elliptic theory it is well-known that the spectrum of ${\frak D}_{M}^{2}$ is discrete and tends to infinity.
We defind the zeta function associated to ${\frak D}_{M}^{2}$ by
$$ \zeta_{{\frak D}_{M}^{2}}(s) = \frac{1}{\Gamma(s)}\int_{0}^{\infty}t^{s-1}Tre^{-t{\frak D}_{M}^{2}} dt,$$
which is holomorphic for $Re(s) > \frac{m}{2}$. 
$\zeta_{{\frak D}_{M}^{2}}(s)$ admits a meromorphic continuation to the whole complex plane having a regular value at $s=0$.
We define the zeta-determinant by
$$Det{\frak D}_{M}^{2} = e^{-\zeta_{{\frak D}_{M}^{2}}^{\prime}(0)}.$$ 

We suppose that $Y$ is a hypersurface of $M$ such that $M-Y$ has two components and $N$ is a collar neighborhood of $Y$ which 
is diffeomorphic to $[-1,1]\times Y$. Choose a metric $g$ on $M$ which is a product metric on $N$.
We now assume the product structures of the bundle $E$ and the Dirac operator ${\frak D}_{M}$ in the following sense.
We first assume that $E|_{N}=p^{\ast}E|_{Y}$, where $p:[-1,1] \times Y \rightarrow Y$ is the canonical projection. 
We also assume that ${\frak D}_{M}$ is of the form on $N$ 
$$ {\frak D}_{M}= G( \partial_{u}+B),$$
where $G : E|_{Y} \rightarrow E|_{Y}$ is a bundle automorphism, $\partial_{u}$ is the normal derivative to $Y$
in the usual direction and $B$ is a Dirac operator on $Y$.
Here, we assume that $G$ and $B$ do not depend on the normal coordinate $u$ and satisfy
$$G^{\ast}=-G, \qquad G^{2}=-Id, \qquad B^{\ast}=B  \qquad {\text and }  \qquad GB=-BG.$$
Then we have, on $N$,
$${\frak D}_{M}^{2}=-\partial_{u}^{2}+B^{2}.$$
We denote by $M_{1}$, $M_{2}$ the closure of each component of $M-Y$ so that
$M_{1}$ ($M_{2}$) contains the part $[-1,0]\times Y$ ($[0,1]\times Y$).
We also denote by ${\frak D}_{M_{1}}$, ${\frak D}_{M_{2}}$ the restriction of ${\frak D}_{M}$ to $M_{1}$, $M_{2}$, respectively.
In this paper we discuss the decomposition of $Det{\frak D}_{M}^{2}$ into contributions coming from $M_{1}$, $M_{2}$ and 
the hypersurface $Y$ in the adiabatic sense. To do this we consider the following adiabatic setting.

\vskip 0.3 true cm

We denote by $M_{r}$ the compact manifold
without boundary obtained by attaching $N_{r+1}=[-r-1,r+1]\times Y$ on 
$M-(-\frac{1}{2},\frac{1}{2})\times Y$ by identifying $[-1,-\frac{1}{2}]\times Y$ with 
$[-r-1,-r-\frac{1}{2}]\times Y$
and $[\frac{1}{2},1]\times Y$ with $[r+\frac{1}{2},r+1]\times Y$. 
We also denote by $M_{1,r}$, $M_{2,r}$ the manifolds with boundary which
are obtained by attaching $[-r,0]\times Y$, $[0,r]\times Y$ on $M_{1}$, $M_{2}$ by identifying $\partial M_{1}$ with $Y_{-r} := \{ -r \}\times Y$
and $\partial M_{2}$ with $Y_{r} := \{ r \}\times Y$, respectively. 
Then the bundle $E\rightarrow M$ and the Dirac operator ${\frak D}_{M}$ on $M$ can be extended naturally to 
the bundle $E_{r}\rightarrow M_{r}$ and the Dirac operator ${\frak D}_{M_{r}}$ on $M_{r}$.  
We also denote by ${\frak D}_{M_{1,r}}$, ${\frak D}_{M_{2,r}}$ ($E_{1,r}$, $E_{2,r}$) the restriction of ${\frak D}_{M_{r}}$
($E_{r}$) to $M_{1,r}$, $M_{2,r}$, respectively
and by  ${\frak D}^{2}_{M_{i,r},D_{0}}$  the Dirac Laplacian ${\frak D}^{2}_{M_{i,r}}$ on $M_{i}$ with the Dirichlet boundary condition on 
$Y_{0} := \{ 0 \} \times Y$,
{\it i.e.}, $Dom({\frak D}^{2}_{M_{i,r},D_{0}})=
\{ \phi \in C^{\infty}(M_{i,r}) \mid \phi|_{Y_{0}}=0 \}$.

We define the operators $Q_{i}:C^{\infty}(Y) \rightarrow C^{\infty}(Y)$ ($i=1$, $2$) as follows. 
For $f\in C^{\infty}(Y)$, choose $\phi_{i}\in C^{\infty}(M_{i})$ satisfying ${\frak D}^{2}_{M_{i}}\phi_{i}=0$ and
$\phi_{i}|_{Y}=f$. Then we define 
$$Q_{1}(f)=(\partial_{u}\phi_{1})|_{Y} , \qquad \qquad Q_{2}(f)=(-\partial_{u}\phi_{2})|_{Y}. \tag1.1 $$
We show in Proposition 4.5 that if both $Q_{1}+\sqrt{B^{2}}$ and $Q_{2}+\sqrt{B^{2}}$ are invertible, then $B$ is invertible for each $r>0$
and  ${\frak D}_{M_{r}}$ is invertible for $r$ large enough. 
The following is the main result of [9] given by the author.

\proclaim{Theorem 1.1}
Let $M$ be a compact oriented Riemannian manifold having the product structures near a hypersurface $Y$.
We assume that both $Q_{1}+\sqrt{B^{2}}$ and $Q_{2}+\sqrt{B^{2}}$ are invertible operators.
Then :
$$
\lim_{r\to\infty}\left\{ \log Det({\frak D}^{2}_{M_{r}})-\log Det({\frak D}^{2}_{M_{1,r},D_{0}})-\log Det({\frak D}^{2}_{M_{2,r},D_{0}}) 
\right\} = \frac{1}{2}\log Det(B^{2}). 
$$
\endproclaim
\noindent
{\it Remark} \hskip 0.1 true cm : \hskip 0.2 true cm
We denote $M_{1,\infty}=M_{1}\cup_{\partial M_{1}}[0,\infty)\times Y$, 
$M_{2,\infty}=M_{2}\cup_{\partial M_{2}}(-\infty,0]\times Y$ and
by ${\frak D}_{M_{i,\infty}}$ ($i=1, 2$) the natural extension of ${\frak D}_{M_{i,r}}$ to $M_{i,\infty}$.
Then we prove in Proposition 4.5 that 
the invertibility of both $Q_{1}+\sqrt{B^{2}}$ and $Q_{2}+\sqrt{B^{2}}$ is equivalent to the non-existence of extended 
$L^{2}$-solutions (see Definition 4.4) of ${\frak D}_{M_{1,\infty}}$ and ${\frak D}_{M_{2,\infty}}$ on 
$M_{1,\infty}$ and $M_{2,\infty}$. \newline

\vskip 0.2 true cm

\noindent
The purpose of this paper is to establish a same type of formula as Theorem 1.1 with 
the Atiyah-Patodi-Singer boundary condition (APS condition) instead of
the Dirichlet boundary condition.

\vskip 0.3 true cm

Recall that $B$ is a Dirac operator on $Y$ and its spectrum is distributed from negative infinity to positive infinity.
Denote by $P_{<}$ ($P_{>}$), $P_{\geq}$ ($P_{\leq}$)  the projections from $C^{\infty}(Y)$ to negative (positive) and non-negative
(non-positive) eigensections of $B$, respectively.
Then ${\frak D}_{M_{1,r},P_{<}}$ and ${\frak D}_{M_{2,r},P_{\geq}}$ are defined by the same operators ${\frak D}_{M_{1,r}}$,
${\frak D}_{M_{2,r}}$ with
$$Dom({\frak D}_{M_{1,r},P_{<}})= \{ \phi\in C^{\infty}(M_{1,r}) \mid P_{<}(\phi|_{Y_{0}})=0 \}, $$
$$Dom({\frak D}_{M_{2,r},P_{\geq}})= \{ \phi\in C^{\infty}(M_{2,r}) \mid P_{\geq}(\phi|_{Y_{0}})=0 \}. $$
Similarly, 
${\frak D}_{M_{1,r},P_{<}}^{2} :=({\frak D}_{M_{1,r},P_{<}})({\frak D}_{M_{1,r},P_{<}})$ and
${\frak D}_{M_{2,r},P_{\geq}}^{2} := ({\frak D}_{M_{2,r},P_{\geq}})({\frak D}_{M_{2,r},P_{\geq}})$
are defined by the same operators ${\frak D}_{M_{1,r}}^{2}$, ${\frak D}_{M_{2,r}}^{2}$ with
$$Dom({\frak D}_{M_{1,r},P_{<}}^{2})= \{ \phi\in C^{\infty}(M_{1,r}) \mid P_{<}(\phi|_{Y_{0}})=0,
P_{\geq}((\partial_{u}\phi + B\phi)|_{Y_{0}})=0 \}, \tag1.2 $$
$$Dom({\frak D}_{M_{2,r},P_{\geq}}^{2})= \{ \phi\in C^{\infty}(M_{2,r}) \mid P_{\geq}(\phi|_{Y_{0}})=0,
P_{<}((\partial_{u}\phi + B\phi)|_{Y_{0}})=0 \}. \tag1.3 $$
\noindent
${\frak D}_{M_{i,r},P_{>}}$, ${\frak D}_{M_{i,r},P_{\leq}}$, ${\frak D}^{2}_{M_{i,r},P_{>}}$ and ${\frak D}^{2}_{M_{i,r},P_{\leq}}$
are defined similarly.

\vskip 0.3 true cm

Put 
$N_{-r,0}=[-r,0]\times Y$ and $N_{0,r}=[0,r]\times Y$. 
Then from the decomposition
$$M_{1,r}=M_{1}\cup_{Y_{-r}} N_{-r,0}, \qquad  M_{2,r}=M_{2}\cup_{Y_{r}} N_{0,r},$$
we have the following theorem, which we call Burghelea-Friedlander-Kappeler's gluing formula 
and refer to [9] for the proof (see also [4], [8]).

\proclaim{Theorem 1.2}
Suppose that $k=dimKer B$. Then :
$$\multline (1) \quad \log Det {\frak D}_{M_{1,r},D_{0}}^{2} = \log Det {\frak D}_{M_{1},D_{-r}}^{2} 
+ \log Det (-\partial_{u}^{2}+B^{2})_{N_{-r,0},D_{-r},D_{0}} \\ 
-\log 2 \cdot (\zeta_{B^{2}}(0)+k)+\log Det R_{M_{1,r},D_{-r}}. \endmultline $$
$$ \multline (2) \quad \log Det {\frak D}_{M_{2,r},D_{0}}^{2} = \log Det {\frak D}_{M_{2},D_{r}}^{2} 
+ \log Det (-\partial_{u}^{2}+B^{2})_{N_{0,r},D_{0},D_{r}}  \\
 -\log 2 \cdot (\zeta_{B^{2}}(0)+k)+\log Det R_{M_{2,r},D_{r}}. \endmultline $$
\endproclaim
\noindent
Here the Dirichlet-to-Neumann operators 
$R_{M_{1,r},D_{-r}} : C^{\infty}(Y_{-r}) \rightarrow C^{\infty}(Y_{-r})$ and 
$R_{M_{2,r},D_{r}} : C^{\infty}(Y_{r}) \rightarrow C^{\infty}(Y_{r})$ are defined as follows.
For $f\in C^{\infty}(Y_{-r})$ and ${\tilde f}\in C^{\infty}(Y_{r})$, choose
$\phi\in C^{\infty}(M_{1})$, $\psi\in C^{\infty}(N_{-r,0})$, 
${\tilde \phi}\in C^{\infty}(M_{2})$, ${\tilde \psi}\in C^{\infty}(N_{0,r})$ so that
$$ 
{\frak D}^{2}_{M_{1}}\phi  = 0, \qquad ( -\partial_{u}^{2}+B^{2}) \psi = 0,  
\qquad \phi|_{Y_{-r}} = \psi|_{Y_{-r}} = f, 
\qquad \psi|_{Y_{0}} = 0 ,$$
$$ 
{\frak D}^{2}_{M_{2}} {\tilde \phi} = 0, \qquad  ( -\partial_{u}^{2}+B^{2}){\tilde \psi} =0, \qquad 
{\tilde \phi}|_{Y_{r}} = {\tilde \psi}|_{Y_{r}} = {\tilde f},
\qquad \quad {\tilde \psi}|_{Y_{0}} = 0 .$$
Then we define :
$$
R_{M_{1,r},D_{-r}}(f) : = (\partial_{u} \phi)|_{Y_{-r}} -  (\partial_{u} \psi)|_{Y_{-r}} = Q_{1}(f) - (\partial_{u} \psi)|_{Y_{-r}} , $$
$$
R_{M_{2,r},D_{r}}({\tilde f}) : = - (\partial_{u} {\tilde \phi})|_{Y_{r}} +  (\partial_{u} {\tilde \psi})|_{Y_{r}} = 
Q_{2}({\tilde f}) + (\partial_{u} {\tilde \psi})|_{Y_{r}} . $$
The operator $(-\partial_{u}^{2}+B^{2})_{N_{-r,0},D_{-r},D_{0}}$ is the Laplacian 
$(-\partial_{u}^{2}+B^{2})$ on $N_{-r,0}$ with the Dirichlet condition on $Y_{-r}$, $Y_{0}$ and 
$(-\partial_{u}^{2}+B^{2})_{N_{0,r},D_{0},D_{r}}$ is defined similarly.

By replacing the Dirichlet condition on $Y_{0}$ by the APS condition, 
we obtain the following theorem, which can be proved by the exactly same way as Theorem 1.2
without any modification.

\proclaim{Theorem 1.3}
Suppose that  ${\frak D}_{M_{1,r},P_{<}}^{2}$ and ${\frak D}_{M_{2,r},P_{\geq}}^{2}$ are invertible operators
and $k=dimKer B$. Then :
$$\multline (1) \quad \log Det {\frak D}_{M_{1,r},P_{<}}^{2} = \log Det {\frak D}_{M_{1},D_{-r}}^{2} 
+ \log Det (-\partial_{u}^{2}+B^{2})_{N_{-r,0},D_{-r},P_{<}} \\ 
-\log 2 \cdot (\zeta_{B^{2}}(0)+k)+\log Det R_{M_{1,r},P_{<}}. \endmultline $$
$$ \multline (2) \quad \log Det {\frak D}_{M_{2,r},P_{\geq}}^{2} = \log Det {\frak D}_{M_{2},D_{r}}^{2} 
+ \log Det (-\partial_{u}^{2}+B^{2})_{N_{0,r},P_{\geq},D_{r}}  \\
 -\log 2 \cdot (\zeta_{B^{2}}(0)+k)+\log Det R_{M_{2,r},P_{\geq}}. \endmultline $$
\endproclaim
\noindent
Here 
$(-\partial_{u}^{2}+B^{2})_{N_{-r,0},D_{-r},P_{<}}$ is the operator $-\partial_{u}^{2}+B^{2}$ on $N_{-r,0}$
with the Dirichlet boundary condition on $Y_{-r}$ and the APS condition $P_{<}$ on $Y_{0}$, {\it i.e.},
$$ \multline
Dom((-\partial_{u}^{2}+B^{2})_{N_{-r,0},D_{-r},P_{<}})=  \\
\{ \phi \in C^{\infty}(N_{-r,0}) \mid  \phi|_{Y_{-r}}=0,
P_{<}(\phi|_{Y_{0}})=0 , P_{\geq}((\partial_{u}\phi+B\phi)|_{Y_{0}})=0  \}, \endmultline $$
and $ (-\partial_{u}^{2}+B^{2})_{N_{0,r},P_{\geq},D_{r}}$ is defined similarly. The operators  \newline
$R_{M_{1,r},P_{<}} : C^{\infty}(Y_{-r}) \rightarrow C^{\infty}(Y_{-r})$ and 
$R_{M_{2,r},P_{\geq}} : C^{\infty}(Y_{r}) \rightarrow C^{\infty}(Y_{r})$, which are the Dirichlet-to-Neumann operators
with the APS condition on $Y_{0}$,
are defined as follows.
For $f \in C^{\infty}(Y_{-r})$ and ${\tilde f}\in C^{\infty}(Y_{r})$, choose $\phi \in C^{\infty}(M_{1})$, $\psi \in C^{\infty}(N_{-r,0})$, 
${\tilde \phi}\in C^{\infty}(M_{2})$, ${\tilde \psi}\in C^{\infty}(N_{0,r})$
satisfying
$$ \multline
{\frak D}^{2}_{M_{1}}\phi=0, \qquad (-\partial_{u}^{2}+B^{2})\psi=0,  \qquad \phi|_{Y_{-r}}=\psi|_{Y_{-r}}=f, \\
P_{<}(\psi|_{Y_{0}})= P_{\geq}((\partial_{u}\psi+B\psi)|_{Y_{0}})=0, \endmultline $$
$$ \multline
{\frak D}^{2}_{M_{2}} {\tilde \phi}=0, \qquad (-\partial_{u}^{2}+B^{2}){\tilde \psi}=0,  \qquad  
{\tilde \phi}|_{Y_{r}} = {\tilde \psi}|_{Y_{r}} = {\tilde f}, \\
P_{\geq}({\tilde \psi}|_{Y_{0}}) =P_{<}((\partial_{u}{\tilde \psi}+B{\tilde \psi})|_{Y_{0}})=0. \endmultline $$
Then we define  
$$R_{M_{1,r},P_{<}}(f) : = (\partial_{u}\phi|)_{Y_{-r}} - (\partial_{u}\psi)|_{Y_{-r}} = Q_{1}(f) - (\partial_{u} \psi)|_{Y_{-r}},$$
$$
R_{M_{2,r},P_{\geq}}({\tilde f}) : = - (\partial_{u} {\tilde \phi})|_{Y_{r}} +  (\partial_{u} {\tilde \psi})|_{Y_{r}} = 
Q_{2}({\tilde f}) + (\partial_{u} {\tilde \psi})|_{Y_{r}} . \quad $$

\vskip 0.3 true cm

Now we discuss a relation  on the cylinder part between the zeta-determinant with the Dirichlet condition and 
the APS condition. With the help of Theorem 1.2 and Theorem 1.3, it is enough to consider the terms coming from
the cylinder parts $[-r,0]\times Y$ and $[0,r]\times Y$.
To describe the main result of this paper, we define the operator
$Q_{(\partial_{u}+|B|),r} : C^{\infty}(Y_{0}) \rightarrow C^{\infty}(Y_{0})$ in the following way, which is suggested in [6].
For $f\in C^{\infty}(Y_{0})$, choose $\phi\in C^{\infty}(N_{-r,0})$ satisfying
$$ (-\partial_{u}^{2}+B^{2}) \phi =0, \qquad \phi|_{Y_{-r}}=0, \qquad \phi|_{Y_{0}}=f.$$
Then we define
$$ Q_{(\partial_{u}+|B|),r} (f) := (\partial_{u}\phi + |B|\phi)|_{Y_{0}}.$$ 
We can check easily by direct computation that  $ Q_{(\partial_{u}+|B|),r}$ is a positive operator (Lemma 4.1, (5)) and
have the following theorem, which is proved in Section 2 and 3 by modifying the proof of Theorem 1.2.
\proclaim{Theorem 1.4}
$$
 \log Det (-\partial_{u}^{2}+B^{2})_{N_{-r,0},D_{-r},P_{<}} + 
\log Det (-\partial_{u}^{2}+B^{2})_{N_{0,r},P_{\geq},D_{r}}  \quad - $$
$$ \log Det (-\partial_{u}^{2}+B^{2})_{N_{-r,0},D_{-r},D_{0}}  - 
\log Det (-\partial_{u}^{2}+B^{2})_{N_{0,r},D_{0},D_{r}}
=  \log Det Q_{(\partial_{u}+|B|),r}.  $$
\endproclaim

Now we are ready to describe the main result of this paper. We prove in Proposition 4.5 that if each $Q_{i}+\sqrt{B^{2}}$
($i=1, 2$) is invertible, then $B$, ${\frak D}^{2}_{M_{1,r},P_{<}}$ and ${\frak D}^{2}_{M_{2,r},P_{>}}$ are invertible for each $r>0$ and
${\frak D}^{2}_{M_{r}}$ is invertible for  $r$ large enough.
Combining Theorem 1.1, 1.2 ,1.3 and 1.4 with the Remark after Theorem 1.1, we have the following theorem,
which is the main result of this paper.

\proclaim{Theorem 1.5}
Let $M$ be a compact oriented Riemannian manifold having the product structures near a hypersurface $Y$.
We assume that there are no extended $L^{2}$-solutions of ${\frak D}_{M_{i,\infty}}$ on $M_{i,\infty}$ for $i=1, 2$.
Then :
$$\lim_{r\to \infty} \{ \log Det {\frak D}_{M_{r}}^{2} - \log Det {\frak D}_{M_{1,r},P_{<}}^{2}
- \log Det {\frak D}_{M_{2,r},P_{>}}^{2} \} 
= - \log 2  \cdot \zeta_{B^{2}}(0) .$$
\endproclaim

\vskip 0.3 true cm

Theorem 1.5 was proved originally in [11] by Park and Wojciechowski on an odd dimensional compact oriented
Riemannian manifold under the same assumption. 
They used the fact that $\zeta_{{\frak D}^{2}_{M_{1},P_{>}}}(0)=0$ when $Ker B=0$ and $dimM$ is odd
and then decomposed $Tr e^{-t{\frak D}^{2}_{M_{r}}}$ into contributions coming from $M_{1,r}$, $M_{2,r}$ and a cylinder part plus
some error terms. They proved that the error terms tend to zero as $r \rightarrow\infty$ and finally
computed the contribution coming from the cylinder part as $r \rightarrow\infty$.
Recently they improved their result in [12] by deleting the assumption of the non-existence of $L^{2}$-solutions on $M_{i,\infty}$ ($i=1,2$). 
For this work they strongly used the scattering theory developed in [10].

In this paper we, however, take different approach for the proof of Theorem 1.5. We are going to use 
Burghelea-Friedlander-Kappeler's gluing formula established in case of the Dirichlet boundary condition.
The original form of this formula contains a constant which can be expressed in terms of zero coefficients of some asymptotic expansions ([4], [8]).
Under the assumption of the product structures near $Y$, it is shown by the author in [9] that this constant is
$-\log 2 \cdot (\zeta_{B^{2}}(0) + dim Ker B)$  and
hence Theorem 1.2 and 1.3 are obtained. We are going to use this result intensively.   
Finally we use Theorem 1.4
to compare the case of the Dirichlet boundary condition with the APS boundary condition
and use Theorem 1.1 to compute the adiabatic limit as $r \rightarrow \infty$. 
One of the advantages for this approach is that this method works in not only odd but also even
dimensional manifolds.

\vskip 1 true cm

\S 2 {\bf Proof of Theorem 1.4 }

\vskip 0.5 true cm

In this section we are going to prove Theorem 1.4. Throughout this section every computation will be done on the cylinders
$[-r,0] \times Y$ and $[0,r]\times Y$.

First of all, one can check by direct computation that the spectra of  \newline
$(-\partial_{u}^{2}+B^{2})_{N_{-r,0},D_{-r},D_{0}}$,
$(-\partial_{u}^{2}+B^{2})_{N_{0,r},D_{0},D_{r}}$,
$(-\partial_{u}^{2}+B^{2})_{N_{-r,0},D_{-r},P_{<}}$ and \newline
$(-\partial_{u}^{2}+B^{2})_{N_{0,r},P_{\geq},D_{r}}$  
 are as follows.
$$(1) \quad Spec\left( (-\partial_{u}^{2}+B^{2})_{N_{-r,0},D_{-r},D_{0}} \right) = 
Spec\left( (-\partial_{u}^{2}+B^{2})_{N_{0,r},D_{0},D_{r}} \right) = 
 \qquad \qquad \qquad \quad $$
$$\{ \lambda^{2}+(\frac{k \pi}{r})^{2} \mid \lambda\in Spec(B),  k = 1,2,3,\cdots \}. $$
$$(2) \quad Spec\left( (-\partial_{u}^{2}+B^{2})_{N_{-r,0},D_{-r},P_{<}} \right)=
\{ \mu_{\lambda, l} \mid \lambda\in Spec(B), \lambda \geq 0, \mu_{\lambda,l}> \lambda^{2} \} \quad \cup 
\qquad \qquad $$
$$ \{ \lambda^{2}+(\frac{k \pi}{r})^{2}  \mid  \lambda\in Spec(B), \lambda < 0, k=1,2,3,\cdots \} . \quad  $$ 
$$(3) \quad Spec\left( (-\partial_{u}^{2}+B^{2})_{N_{0,r},P_{\geq},D_{r}} \right)=
\{ \mu_{\lambda, l} \mid \lambda\in Spec(B), \lambda > 0, \mu_{\lambda,l}> \lambda^{2} \} \quad \cup
\qquad \qquad $$
$$ \{ \lambda^{2}+(\frac{k \pi}{r})^{2} \mid \lambda\in Spec(B), \lambda \geq 0, k=1,2,3,\cdots \} .$$
In (2) and (3) $\mu_{\lambda, l}$'s are the solutions of the equation
$$\sqrt{\mu-\lambda^{2}} \hskip 0.1 true cm cos(\sqrt{\mu-\lambda^{2}} \hskip 0.1 true cm r)
+\lambda \hskip 0.1 true cm sin(\sqrt{\mu-\lambda^{2}} \hskip 0.1 true cm r) = 0.$$

We next introduce the boundary condition $\partial_{u}+|B|$ on $Y_{0}$ and consider the Laplacian
$ (-\partial_{u}^{2}+B^{2})_{N_{-r,0},D_{-r},\partial_{u}+|B|}$ with 
$$ \multline
 Dom\left((-\partial_{u}^{2}+B^{2})_{N_{-r,0},D_{-r},\partial_{u}+|B|}\right) =  \\ 
\{ \phi\in C^{\infty}(N_{-r,0}) \mid
\phi|_{Y_{-r}}=0, (\partial_{u}\phi +|B|\phi)|_{Y_{0}}=0 \}. \endmultline $$
Then the spectrum of this operator is :
$$ (4) \quad 
 Spec\left((-\partial_{u}^{2}+B^{2})_{N_{-r,0},D_{-r},\partial_{u}+|B|}\right) =  
 \{ \mu_{\lambda, l} \mid   \lambda\in Spec(|B|),  
\mu_{\lambda,l}> \lambda^{2} \}, $$
where $\mu_{\lambda,l}$'s are the solutions of the equation
$$ \sqrt{\mu-\lambda^{2}} \hskip 0.1 true cm cos(\sqrt{\mu-\lambda^{2}} \hskip 0.1 true cm r)
+\lambda \hskip 0.1 true cm sin(\sqrt{\mu-\lambda^{2}} \hskip 0.1 true cm r) = 0.$$

\vskip 0.2 true cm

Hence from (1), (2), (3), (4), we have
$$\multline
\zeta_{(-\partial_{u}^{2}+B^{2})_{N_{-r,0},D_{-r},P_{<}}}(s)+
\zeta_{(-\partial_{u}^{2}+B^{2})_{N_{0,r},P_{\geq},D_{r}}}(s) \\
-\zeta_{(-\partial_{u}^{2}+B^{2})_{N_{-r,0},D_{-r},D_{0}}}(s)-
\zeta_{(-\partial_{u}^{2}+B^{2})_{N_{0,r},D_{0},D_{r}}}(s) \endmultline  $$
$$ 
= \hskip 0.1 true cm \zeta_{(-\partial_{u}^{2}+B^{2})_{N_{-r,0},D_{-r},(\partial_{u}+|B|)}}(s) - 
\zeta_{(-\partial_{u}^{2}+B^{2})_{N_{-r,0},D_{-r},D_{0}}}(s). \tag2.1  $$

\vskip 0.3 true cm

Now we are going to use the method in [6] and [8] to analyze
$$ \log Det( (-\partial_{u}^{2}+B^{2})_{N_{-r,0},D_{-r},(\partial_{u}+|B|)} ) - 
\log Det( (-\partial_{u}^{2}+B^{2})_{N_{-r,0},D_{-r},D_{0}} ).$$
From now on, we simply denote $(-\partial_{u}^{2}+B^{2})$, $\partial_{u}+|B|$ by 
$\Delta$, $C$ and the bundle $E|_{N_{-r,0}}$ by $E$.
We first consider $\Delta^{l}+t^{l}$ for any positive integer $l>[\frac{m}{2}]$ and $t>0$ rather than $\Delta$ itself
because under some proper boundary condition $(\Delta^{l}+t^{l})^{-1}$ is a trace class operator and in this case we
are able to use the well-known formula about the derivative of $\log Det(\Delta^{l}+t^{l})$ with respect to $t$ ({\it c.f.} (2.3)).
For simplicity we put $l=m$. 
Note that 
$$
\Delta^{m}+t^{m}=\cases \prod_{k=-[\frac{m}{2}]}^{[\frac{m-1}{2}]}
(\Delta+e^{i\frac{2k\pi}{m}}t),
& \text{ if $m$ is odd} \\ 
\prod_{k=-[\frac{m}{2}]}^{[\frac{m-1}{2}]}(\Delta+e^{i\frac{(2k+1)\pi}{m}}t),
& \text{ if $m$ is even}.
\endcases
$$
For $-[\frac{m}{2}]\leq k \leq [\frac{m-1}{2}]$, denote
$$
\alpha_{0}=\cases e^{-i\frac{2\pi}{m}[\frac{m}{2}]}, & \text{ if $m$ is odd} \\ 
e^{-i\frac{(m-1)\pi}{m}}, & \text{ if $m$ is even}
\endcases
$$
and $\alpha_{k}=\alpha_{0}e^{i\frac{2k\pi}{m}}$ ($k=0,1,2,\cdots, m-1$).

Suppose that $\gamma_{-r}$ and $\gamma_{0}$ are the restriction operators from $C^{\infty}(E)$ to $C^{\infty}(Y_{-r})$
and $C^{\infty}(Y_{0})$, respectively. Then the Dirichlet boundary conditions $D_{-r}$ and $D_{0}$ on $Y_{-r}$ and $Y_{0}$
are defined by the operators $\gamma_{-r}$ and $\gamma_{0}$.
We can check easily that for each $\alpha_{k}$ and $t>0$,
$$
(\Delta+\alpha_{k}t)_{D_{-r},D_{0}} : \{ \phi\in C^{\infty}(E) \mid \gamma_{-r}\phi = \gamma_{0}\phi =0 \} \rightarrow C^{\infty}(E) $$
is an invertible operator and we can define the Poisson operator 
$P_{D_{0}}(\alpha_{k}t) : C^{\infty}(E|_{Y_{0}}) \rightarrow C^{\infty}(E)$, which is characterized as follows.
$$\gamma_{-r} P_{D_{0}}(\alpha_{k}t) =0, \qquad \gamma_{0}P_{D_{0}}(\alpha_{k}t) = Id_{Y_{0}}, \qquad 
(\Delta + \alpha_{k}t)P_{D_{0}}(\alpha_{k}t) =0. $$ 

Now we are going to define boundary conditions corresponding to the operator $\Delta^{m}+t^{m}$.
Define
$D_{-r,m}(t)$, $D_{0,m}(t)$ and $C_{m}(t) : C^{\infty}(E) \rightarrow \oplus_{m}C^{\infty}(E|_{Y_{0}}) $
as follows.
$$ 
D_{-r,m}(t)= $$
$$\left(\gamma_{-r}, \gamma_{-r}(\Delta+\alpha_{0}t), \gamma_{-r} (\Delta+\alpha_{1}t) (\Delta+\alpha_{0}t),
\cdots , \gamma_{-r}(\Delta+\alpha_{m-2}t)  \cdots (\Delta+\alpha_{0}t) \right),
$$
$$ 
D_{0,m}(t)= $$
$$\left(\gamma_{0}, \gamma_{0}(\Delta+\alpha_{0}t), \gamma_{0} (\Delta+\alpha_{1}t) (\Delta+\alpha_{0}t),
\cdots , \gamma_{0}(\Delta+\alpha_{m-2}t)  \cdots (\Delta+\alpha_{0}t) \right),
$$
$$ 
C_{m}(t)= $$
$$\left(\gamma_{0}C, \gamma_{0}C(\Delta+\alpha_{0}t), \gamma_{0}C (\Delta+\alpha_{1}t) (\Delta+\alpha_{0}t),
\cdots , \gamma_{0}C(\Delta+\alpha_{m-2}t)  \cdots (\Delta+\alpha_{0}t) \right).
$$
Then the Poisson operator 
${\tilde P_{D_{0,m}(t)}}(t) : \oplus_{m}C^{\infty}(E|_{Y_{0}})\rightarrow C^{\infty}(E)$ associated to 
$(\Delta^{m}+t^{m},D_{0,m}(t))$ is given as follows ({\it c.f.} [4], [8]).
$$
\split
{\tilde P_{D_{0,m}(t)}}(t)(f_{0}, \cdots, f_{m-1})=P_{D_{0}}(\alpha_{0}t)f_{0}+
(\Delta+\alpha_{0}t)_{D_{-r},D_{0}}^{-1}P_{D_{0}}(\alpha_{1}t)f_{1} + \cdots \\
+ (\Delta+\alpha_{0}t)_{D_{-r},D_{0}}^{-1}(\Delta+\alpha_{1}t)_{D_{-r},D_{0}}^{-1}\cdots
(\Delta+\alpha_{m-2}t)_{D_{-r},D_{0}}^{-1}P_{D_{0}}(\alpha_{m-1}t)f_{m-1}.
\endsplit
$$
{\it I.e.}  ${\tilde P_{D_{0,m}(t)}}(t)$ defined as above satisfies the following properties.
$$D_{-r,m}(t) {\tilde P_{D_{0,m}(t)}}(t) = 0 ,  \qquad   D_{0,m}(t) {\tilde P_{D_{0,m}(t)}}(t) = Id_{\oplus_{m}C^{\infty}(E|_{Y_{0}})} , $$
$$ \text{ and } \quad (\Delta^{m}+t^{m}) {\tilde P_{D_{0,m}(t)}}(t) = 0.$$
Note that for $i=-r$, $0$,
$$ D_{i,m}(0)= (\gamma_{i}, \gamma_{i}\Delta, \cdots, \gamma_{i}\Delta^{m-1}), \qquad
 C_{m}(0)= (\gamma_{0}C, \gamma_{0}C\Delta, \cdots, \gamma_{0}C\Delta^{m-1}).$$
Put
$$ \Omega(t) = $$

$$
\left(\smallmatrix
1, & 0, & 0, & \cdots ,& 0  , & 0 \\
-\alpha_{0}t, & 1, & 0, & \cdots , & 0 , & 0 \\
\alpha_{0}^{2}t^{2}, & -t(\alpha_{0} + \alpha_{1}), & 1, & \cdots ,& 0 , & 0 \\
\vdots & \vdots & \vdots & \vdots & \vdots & \vdots \\
(-1)^{m-1}\alpha_{0}^{m-1}t^{m-1}, \quad &
(-1)^{m-2}t^{m-2}\sum_{k=0}^{m-2}\alpha_{0}^{m-2-k}\alpha_{1}^{k}, \quad &
\cdots , \qquad & \cdots , \quad & -t \sum_{i=0}^{m-2} \alpha_{i} , \quad  & 1
\endsmallmatrix \right)
$$
Then one can check by direct computation that 
$$ D_{-r,m}(0)=\Omega(t)D_{-r,m}(t), \quad     
D_{0,m}(0)=\Omega(t)D_{0,m}(t) \quad \text{ and } \quad C_{m}(0)=\Omega(t)C_{m}(t).$$
Define $P_{D_{0,m}(0)}(t)={\tilde P}_{D_{0,m}(t)}(t)\Omega(t)^{-1}$. Then
$P_{D_{0,m}(0)}(t)$ is the Poisson operator associated to
$(\Delta^{m}+t^{m}, D_{0,m}(0))$, which is characterized as follows.
$$
D_{-r,m}(0)P_{D_{0,m}(0)}(t)=0, \qquad  D_{0,m}(0)P_{D_{0,m}(0)}(t)= Id_{\oplus_{m} C^{\infty}(Y_{0})}, $$
$$ \text{ and }  (\Delta^{m}+t^{m})P_{D_{0,m}(0)}(t)=0.$$

Next we consider 
$$(\Delta^{m}+t^{m})_{D_{-r,m}(0),D_{0,m}(0)} : 
\{ \phi\in C^{\infty}(E) \mid D_{-r,m}(0) \phi = D_{0,m}(0) \phi = 0 \} \rightarrow C^{\infty}(E)$$
and 
$$(\Delta^{m}+t^{m})_{D_{-r,m}(0),C_{m}(0)} : 
\{ \phi\in C^{\infty}(E) \mid D_{-r,m}(0) \phi = C_{m}(0) \phi = 0 \} \rightarrow C^{\infty}(E),$$
both of which are invertible operators and $(\Delta^{m}+t^{m})_{D_{-r,m}(0),D_{0,m}(0)}^{-1}$ and \newline 
$(\Delta^{m}+t^{m})_{D_{-r,m}(0),C_{m}(0)}^{-1}$ are trace class operators.
From the following identities

$$
(\Delta^{m}+t^{m})\lbrace (\Delta^{m}+t^{m})_{D_{-r,m}(0),C_{m}(0)}^{-1} - (\Delta^{m}+t^{m})_{D_{-r,m}(0),D_{0,m}(0)}^{-1} \rbrace
= 0,$$

$$D_{-r,m}(0) \lbrace (\Delta^{m}+t^{m})_{D_{-r,m}(0),C_{m}(0)}^{-1} - (\Delta^{m}+t^{m})_{D_{-r,m}(0),D_{0,m}(0)}^{-1} \rbrace
= 0,$$

$$ \multline
D_{0,m}(0) \lbrace (\Delta^{m}+t^{m})_{D_{-r,m}(0),C_{m}(0)}^{-1} - (\Delta^{m}+t^{m})_{D_{-r,m}(0),D_{0,m}(0)}^{-1} \rbrace \\
= D_{0,m}(0)  (\Delta^{m}+t^{m})_{D_{-r,m}(0),C_{m}(0)}^{-1},
\endmultline $$
\newline
we have 

$$ \multline
(\Delta^{m}+t^{m})_{D_{-r,m}(0),C_{m}(0)}^{-1} - (\Delta^{m}+t^{m})_{D_{-r,m}(0),D_{0,m}(0)}^{-1} = \\
 P_{D_{0,m}(0)}(t) D_{0,m}(0) (\Delta^{m}+t^{m})_{D_{-r,m}(0),C_{m}(0)}^{-1}. \endmultline \tag2.2 $$
\noindent
Hence, we have

$$
\frac{d}{dt} \lbrace \log Det (\Delta^{m}+t^{m})_{D_{-r,m}(0),C_{m}(0)} - 
\log Det (\Delta^{m}+t^{m})_{D_{-r,m}(0),D_{0,m}(0)} \rbrace $$
\vskip -2mm
$$ = Tr \lbrace mt^{m-1} \left( (\Delta^{m}+t^{m})_{D_{-r,m}(0),C_{m}(0)}^{-1} - 
(\Delta^{m}+t^{m})_{D_{-r,m}(0),D_{0,m}(0)}^{-1} \right) \rbrace $$
\vskip -2mm
$$= mt^{m-1} Tr \lbrace 
P_{D_{0,m}(0)}(t) D_{0,m}(0) (\Delta^{m}+t^{m})_{D_{-r,m}(0),C_{m}(0)}^{-1} \rbrace 
\qquad \qquad \qquad  $$
\vskip -2mm
$$ = mt^{m-1} Tr \lbrace 
 D_{0,m}(0) (\Delta^{m}+t^{m})_{D_{-r,m}(0),C_{m}(0)}^{-1}P_{D_{0,m}(0)}(t) \rbrace 
\qquad \qquad  \tag2.3 $$

\vskip 0.2 true cm

We now define $Q_{(\partial_{u}+|B|),r} : C^{\infty}(Y_{0}) \rightarrow C^{\infty}(Y_{0})$ 
as in Section 1 by

$$  Q_{(\partial_{u}+|B|),r} = \gamma_{0} (\partial_{u}+|B|) P_{D_{0}}(0).$$ 
\noindent
and for simplicity we denote $Q_{(\partial_{u}+|B|),r}$ by $Q$.
We also define $Q_{m}(t)$, ${\tilde Q}_{m}(t) : \oplus_{m} C^{\infty}(Y_{0}) \rightarrow \oplus_{m} C^{\infty}(Y_{0}) $ as follows.

$$Q_{m}(t)=C_{m}(0) P_{D_{0,m}(0)}(t), \qquad {\tilde Q_{m}}(t)=C_{m}(t) {\tilde P_{D_{0,m}(t)}}(t).$$
\noindent
Then

$$ \split
{\tilde Q_{m}}(t) & = \Omega(t)^{-1} C_{m}(0) P_{D_{0,m}(0)}(t) \Omega(t) \\
& = \Omega(t)^{-1} Q_{m}(t) \Omega(t),
\endsplit \tag2.4 $$
\noindent
and hence $Q_{m}(t)$ and ${\tilde Q}_{m}(t)$ are isospectral and have the same determinants.

Now we are going to describe $ \frac{d}{dt} Q_{m}(t)=C_{m}(0) \frac{d}{dt}P_{D_{0,m}(0)}(t).$
\proclaim{Lemma 2.1}
$$\frac{d}{dt}P_{D_{0,m}(0)}(t) = -m t^{m-1} (\Delta^{m}+t^{m})_{D_{-r,m}(0),D_{0,m}(0)}^{-1} P_{D_{0,m}(0)}(t).$$
\endproclaim
\noindent {\it Proof} \hskip 0.3 true cm
Differentiating $(\Delta^{m}+t^{m}) P_{D_{0,m}(0)}(t) = 0$, we have
$$ (\Delta^{m}+t^{m}) \frac{d}{dt} P_{D_{0,m}(0)}(t) = -m t^{m-1} P_{D_{0,m}(0)}(t). \tag2.5$$
From the following identities
$$D_{-r,m}(0) P_{D_{0,m}(0)}(t) = 0, \qquad D_{0,m}(0) P_{D_{0,m}(0)}(t) = Id_{\oplus_{m}C^{\infty}(Y_{0})},$$
we have
$$
D_{-r,m}(0) \frac{d}{dt}P_{D_{0,m}(0)}(t) = 0, \qquad  D_{0,m}(0) \frac{d}{dt} P_{D_{0,m}(0)}(t) = 0. \tag2.6$$
From (2.5) and (2.6), the result follows.
\qed

From Lemma 2.1 and (2.2),
$$ 
\frac{d}{dt}Q_{m}(t) = -m t^{m-1} C_{m}(0) 
(\Delta^{m}+t^{m})_{D_{-r,m}(0),D_{0,m}(0)}^{-1} P_{D_{0,m}(0)}(t)$$
$$=m t^{m-1} C_{m}(0) \left( (\Delta^{m}+t^{m})_{D_{-r,m}(0),C_{m}(0)}^{-1} - (\Delta^{m}+t^{m})_{D_{-r,m}(0),D_{0,m}(0)}^{-1} \right)
P_{D_{0,m}(0)}(t) $$
$$=m t^{m-1} C_{m}(0) P_{D_{0,m}(0)}(t) D_{0,m}(0) (\Delta^{m}+t^{m})_{D_{-r,m}(0),C_{m}(0)}^{-1} P_{D_{0,m}(0)}(t) \qquad \qquad \qquad $$
$$=m t^{m-1} Q_{m}(t) D_{0,m}(0) (\Delta^{m}+t^{m})_{D_{-r,m}(0),C_{m}(0)}^{-1} P_{D_{0,m}(0)}(t). 
\qquad \qquad \qquad \qquad \qquad $$
As a consequence, we have
$$ \frac{d}{dt}\log Det Q_{m}(t) = Tr\left( Q_{m}(t)^{-1} \frac{d}{dt}Q_{m}(t) \right) $$
$$ =  m t^{m-1} Tr\left( D_{0,m}(0) (\Delta^{m}+t^{m})_{D_{-r,m}(0),C_{m}(0)}^{-1} P_{D_{0,m}(0)}(t)\right). 
 \tag2.7 $$
From (2.3), (2.4) and (2.7), we have
$$   
\frac{d}{dt}\log Det {\tilde Q}_{m}(t) = $$
$$ \frac{d}{dt} \left\{ \log Det (\Delta^{m}+t^{m})_{D_{-r,m}(0),C_{m}(0)} -
\log Det (\Delta^{m}+t^{m})_{D_{-r,m}(0),D_{0,m}(0)} \right\}.  \tag2.8 $$
Setting $Q_{(\partial_{u}+|B|),r}(\alpha_{k}t)=\gamma_{0} C P_{D_{0}}(\alpha_{k}t) : 
C^{\infty}(E|_{Y_{0}})\rightarrow C^{\infty}(E|_{Y_{0}})$ (briefly $Q(\alpha_{k}t)$),
${\tilde Q_{m}}(t)$ is of the following upper triangular matrix,

$$
\left( \smallmatrix
Q(\alpha_{0}t), & \gamma_{0} C (\Delta+\alpha_{0}t)^{-1}_{D_{-r},D_{0}}P_{D_{0}}(\alpha_{1}t),
& \hdots, & \gamma_{0} C (\Delta+\alpha_{0}t)^{-1}_{D_{-r},D_{0}}\hdots 
(\Delta+\alpha_{m-2}t)^{-1}_{D_{-r},D_{0}}P_{D_{0}}(\alpha_{m-1}t) \\
0, & Q(\alpha_{1}t), & \hdots, 
& \gamma_{0} C (\Delta+\alpha_{1}t)^{-1}_{D_{-r},D_{0}}
\hdots (\Delta+\alpha_{m-2}t)^{-1}_{D_{-r},D_{0}}P_{D_{0}}(\alpha_{m-1}t) \\
\vdots & \vdots & \ddots & \vdots \\
0, & 0, & \hdots, & Q(\alpha_{m-1}t) 
\endsmallmatrix \right)
$$
\noindent
and hence we have
$$ 
\frac{d}{dt}\left( \log Det {\tilde Q_{m}}(t) \right)
=\frac{d}{dt}(\sum_{k=0}^{m-1}\log DetQ(\alpha_{k}t)).   \tag2.9 $$
Finally by (2.8), (2.9), we have
$$ \multline 
\log Det (\Delta^{m}+t^{m})_{D_{-r,m}(0),C_{m}(0)} -
\log Det (\Delta^{m}+t^{m})_{D_{-r,m}(0),D_{0,m}(0)} =  \\ 
{\tilde c} + \sum_{k=0}^{m-1} \log Det Q(\alpha_{k}t),
\endmultline \tag2.10 $$
where ${\tilde c}$ does not depend on $t$.
It is known in [4] that $\log Det (\Delta^{m}+t^{m})_{D_{-r,m}(0),C_{m}(0)}$, 
$\log Det (\Delta^{m}+t^{m})_{D_{-r,m}(0),D_{0,m}(0)}$
and $\log DetQ(\alpha_{k}t)$ have asymptotic expansions as $t\rightarrow \infty$ 
and the zero coefficients in the asymptotic expansions
of $\log Det (\Delta^{m}+t^{m})_{D_{-r,m}(0),C_{m}(0)}$ and 
$\log Det (\Delta^{m}+t^{m})_{D_{-r,m}(0),D_{0,m}(0)}$ are zeros ([8], [13]).
Denoting by $c_{k}$ the zero coefficient in the asymptotic expansion of $\log Det  Q(\alpha_{k}t)$,
$${\tilde c}=-\sum_{k=0}^{m-1} c_{k}.$$
Setting $t=0$, we have the following theorem.
\proclaim{Theorem 2.2}
$$\log Det \Delta_{D_{-r},C} - \log Det \Delta_{D_{-r},D_{0}}
= - c + \log Det Q_{(\partial_{u}+|B|),r} ,$$
where $ c = \frac{1}{m}\sum_{k=0}^{m-1} c_{k}$.
\endproclaim
\noindent
{\it Proof} \hskip 0.3 true cm
Setting $t=0$ in (2.10), we have
$$ \multline
\log Det (\Delta^{m})_{D_{-r,m}(0),C_{m}(0)} -
\log Det (\Delta^{m})_{D_{-r,m}(0),D_{0,m}(0)} = \\   
{\tilde c} + m \log Det Q_{(\partial_{u}+|B|),r} . \endmultline $$
From the fact that $\lambda$ is an eigenvalue of $\Delta_{D_{-r},D_{0}}$ ($\Delta_{D_{-r},C}$)
if and only if $\lambda^{m}$ is an eigenvalue of $(\Delta^{m})_{D_{-r,m}(0),D_{0,m}(0)}$
($(\Delta^{m})_{D_{-r,m}(0),C_{m}(0)}$), the result follows.
\qed \newline
From (2.1) and Theorem 2.2, we have the following result.
\proclaim{Theorem 2.3}
$$ \multline
\log Det(-\partial_{u}^{2}+B^{2})_{N_{-r,0},D_{-r},P_{<}}+
\log Det(-\partial_{u}^{2}+B^{2})_{N_{0,r},P_{\geq},D_{r}} \quad - \\
\log Det(-\partial_{u}^{2}+B^{2})_{N_{-r,0},D_{-r},D_{0}}- 
\log Det(-\partial_{u}^{2}+B^{2})_{N_{0,r},D_{0},D_{r}}  \endmultline $$
$$ = - c + \log Det Q_{(\partial_{u}+|B|),r}.  \qquad \qquad \qquad \qquad \qquad \qquad \qquad \qquad \qquad \qquad \qquad $$
\endproclaim  

In the next section, we are going to show that this constant $c=0$ by using the method in [13], which completes the proof of Theorem 1.4

\vskip 1 true cm

\S 3 {\bf Computation of the constant term in Theorem 2.3 }

\vskip 0.5 true cm

Recall that $c=\frac{1}{m}\sum_{k=0}^{m-1}c_{k}$, where $c_{k}$ is the zero coefficient in the asymptotic expansion of 
$\log Det Q_{(\partial_{u}+|B|),r}(\alpha_{k}t) $ as $t\rightarrow\infty$.
Note that for $f\in C^{\infty}(Y)$ with $Bf=\lambda f$,
$$ 
Q_{(\partial_{u}+|B|),r}(\alpha_{k}t) f= \left( \sqrt{\lambda^{2}+\alpha_{k}t \hskip 0.1 true cm}+|\lambda|
+\frac{2 \hskip 0.1 true cm \sqrt{\lambda^{2}+\alpha_{k}t \hskip 0.1 true cm} \hskip 0.1 true cm 
e^{-\sqrt{\lambda^{2}+\alpha_{k}t \hskip 0.1 true cm} \hskip 0.1 true cm r}}
{e^{\sqrt{\lambda^{2}+\alpha_{k}t \hskip 0.1 true cm} \hskip 0.1 true cm r} -
e^{-\sqrt{\lambda^{2}+\alpha_{k}t \hskip 0.1 true cm} \hskip 0.1 true cm r}}\right)f. $$
Here we take the negative real axis as a branch cut for square root of $\alpha_{k}$.
Since the asymptotic expansion of $\log Det Q_{(\partial_{u}+|B|),r}(\alpha_{k}t) $ as $t\rightarrow\infty$
is completely determined up to smoothing operators 
(c.f. [4]), $\log Det Q_{(\partial_{u}+|B|),r}(\alpha_{k}t) $ and  \newline
$\log Det (\sqrt{B^{2}+\alpha_{k}t} +|B|)$ have the same asymptotic expansions.
Hence it's enough to consider the asymptotic expansion of  $\log Det (\sqrt{B^{2}+\alpha_{k}t} +|B|)$.

\vskip 0.3 true cm

Since $Re(\alpha_{k})$ can be negative, to avoid this difficulty we 
choose an angle $\phi_{k}$ so that $0\leq |\phi_{k}| < \frac{\pi}{2}$ and $Re(e^{i(\theta_{k}-\phi_{k})})>0$, where
$\alpha_{k}=e^{i\theta_{k}}$. Then,
$$ \log Det (\sqrt{B^{2}+\alpha_{k}t \hskip 0.1 true cm} +|B|) \qquad \qquad \qquad \qquad \qquad $$
$$= \log Det \lbrace e^{i\frac{\phi_{k}}{2}} 
(\sqrt{e^{-i\phi_{k}}B^{2}}+ \sqrt{e^{-i\phi_{k}}B^{2}+e^{i(\theta_{k}-\phi_{k})}t \hskip 0.1 true cm} \hskip 0.1 true cm ) \rbrace 
\qquad \qquad \qquad \qquad  \quad $$
$$=-\frac{d}{ds}|_{s=0} \lbrace e^{-i\frac{\phi_{k}}{2}s}
\zeta_{\left(\sqrt{e^{-i\phi_{k}}B^{2}}+ \sqrt{e^{-i\phi_{k}}B^{2}+e^{i(\theta_{k}-\phi_{k})}t \hskip 0.1 true cm}
\hskip 0.1 true cm \right) }(s) \rbrace
\qquad \qquad \qquad \qquad \quad $$
$$ \multline
=\frac{i\phi_{k}}{2}\zeta_{\left(\sqrt{e^{-i\phi_{k}}B^{2}}+ 
\sqrt{e^{-i\phi_{k}}B^{2}+e^{i(\theta_{k}-\phi_{k})}t \hskip 0.1 true cm} \hskip 0.1 true cm \right) }(0) \quad + \\
\log Det \left(\sqrt{e^{-i\phi_{k}}B^{2}}+ \sqrt{e^{-i\phi_{k}}B^{2}+e^{i(\theta_{k}-\phi_{k})}t \hskip 0.1 true cm}
\hskip 0.1 true cm \right). \endmultline 
\tag3.1 $$

Now put ${\tilde \theta_{k}}=\theta_{k}-\phi_{k}$ and denote 
$\zeta_{\left(\sqrt{e^{-i\phi_{k}}B^{2}}+ 
\sqrt{e^{-i\phi_{k}}B^{2}+e^{i(\theta_{k}-\phi_{k})}t \hskip 0.1 true cm} \hskip 0.1 true cm \right) }(s)$ simply by $\zeta(s)$.
Then,
$$
\zeta(s)=\frac{1}{\Gamma(s)}\int_{0}^{\infty}r^{s-1}
Tr e^{-r\left(\sqrt{e^{-i\phi_{k}}B^{2}}+ 
\sqrt{e^{-i\phi_{k}}B^{2}+e^{i{\tilde \theta_{k}}}t \hskip 0.1 true cm} \hskip 0.1 true cm \right) } dr $$
$$=\sum_{q=0}^{\infty}\frac{(-1)^{q}}{q!}\frac{1}{\Gamma(s)}\int_{0}^{\infty}r^{s+q-1}
Tr\lbrace \sqrt{e^{-i\phi_{k}}B^{2}}^{q} 
e^{-r\sqrt{e^{-i\phi_{k}}B^{2}+e^{i{\tilde \theta_{k}}}t \hskip 0.1 true cm} \hskip 0.1 true cm} \rbrace dr .$$
Setting 
$$\zeta_{q}(s)= \frac{1}{\Gamma(s)}\int_{0}^{\infty}r^{s+q-1}
Tr\lbrace \sqrt{e^{-i\phi_{k}}B^{2}}^{q} e^{-r\sqrt{e^{-i\phi_{k}}B^{2}+
e^{i{\tilde \theta_{k}}}t \hskip 0.1 true cm} \hskip 0.1 true cm} \rbrace dr , \tag3.2 $$
$$
\zeta_{q}(s) = \frac{1}{\Gamma(s)} \sum_{\lambda_{j}\in Spec(|B|)}
\frac{\sqrt{e^{-i\phi_{k}}\lambda_{j}^{2} \hskip 0.1 true cm}^{q}}{z_{j}^{s+q}} \int_{0}^{\infty}(rz_{j})^{s+q-1} 
e^{-rz_{j}}z_{j} dr, $$
where $z_{j}= \sqrt{e^{-i\phi_{k}}\lambda_{j}^{2}+e^{i{\tilde \theta_{k}}}t}.$
Consider the contour integral $\int_{C}z^{s+q-1}e^{-z}dz$ for $Re s > -q$, where for $arg(z_{j})=\rho_{j}$,
$$C= \{ re^{i\rho_{j}} \mid \epsilon \leq r\leq R \}  \cup  \{ \epsilon e^{i\theta} \mid 0\leq \theta\leq \rho_{j} \}
\cup \{ r \mid \epsilon \leq r\leq R \} \cup  \{ Re^{i\theta} \mid 0\leq \theta\leq \rho_{j} \}$$
and oriented counterclockwise.
Then one can check that
$$ \int_{0}^{\infty}(rz_{j})^{s+q-1} e^{-rz_{j}}z_{j} dr = \Gamma(s+q),  \tag3.3 $$
and hence we have
$$\zeta_{q}(s)= \frac{\Gamma(s+q)}{\Gamma(s)}\sum_{\lambda_{j}\in Spec(|B|)}
\frac{\sqrt{e^{-i\phi_{k}}\lambda_{j}^{2} \hskip 0.1 true cm}^{q}}{\sqrt{e^{-i\phi_{k}}\lambda_{j}^{2}+
e^{i{\tilde \theta_{k}}}t \hskip 0.1 true cm}^{s+q}}. \qquad \qquad \qquad \qquad   $$
Using the equation (3.3) again, we obtain that
$$\zeta_{q}(s)=\frac{\Gamma(s+q)}{\Gamma(s)} \frac{1}{\Gamma(\frac{s+q}{2})}
\int_{0}^{\infty}r^{\frac{s+q}{2}-1}
Tr\lbrace \sqrt{e^{-i\phi_{k}}B^{2}}^{q} e^{-r(e^{-i\phi_{k}}B^{2}+e^{i{\tilde \theta_{k}}}t)} \rbrace dr. $$
Putting $rt=u$,
$$
\zeta_{q}(s)= \frac{\Gamma(s+q)}{\Gamma(\frac{s+q}{2})\Gamma(s)} \hskip 0.1 true cm t^{-\frac{s+q}{2}}
\int_{0}^{\infty} u^{\frac{s+q}{2}-1} e^{-ue^{i{\tilde \theta_{k}}}}
Tr\lbrace \sqrt{e^{-i\phi_{k}}B^{2}}^{q} e^{-\frac{u}{t}e^{-i\phi_{k}}B^{2}} \rbrace du.
\tag3.4 $$
It is known in [7] (see also [3]) that as $r\rightarrow 0$
$$Tr \lbrace \sqrt{e^{-i\phi_{k}}B^{2}}^{q} e^{-re^{-i\phi_{k}}B^{2}} \rbrace \sim
\sum_{j=0}^{\infty}a_{j}r^{\frac{j-(m-1)-q}{2}} + \sum_{j=0}^{\infty}(b_{j}\log r + c_{j})r^{j}.\tag3.5 $$
Hence, as $t\rightarrow\infty$,
$$
\zeta_{q}(s) \sim \frac{\Gamma(s+q)}{\Gamma(\frac{s+q}{2})\Gamma(s)} \hskip 0.1 true cm \lbrace
\hskip 0.1 true cm \sum_{j=0}^{\infty}a_{j}\int_{0}^{\infty}u^{\frac{s+j-(m-1)}{2}-1}e^{-ue^{i{\tilde \theta_{k}}}}
 du \cdot t^{\frac{(m-1)-j-s}{2}} $$
$$+ \hskip 0.1 true cm \sum_{j=0}^{\infty}b_{j}\int_{0}^{\infty}u^{\frac{s+q}{2}+j-1}e^{-ue^{i{\tilde \theta_{k}}}}
\log u  \hskip 0.1 true cm du \cdot t^{-\frac{s+q}{2}-j} $$
$$- \hskip 0.1 true cm \sum_{j=0}^{\infty}b_{j}\int_{0}^{\infty}u^{\frac{s+q}{2}+j-1}e^{-ue^{i{\tilde \theta_{k}}}}
du \cdot t^{-\frac{s+q}{2}-j}\log t \quad $$
$$+ \hskip 0.1 true cm \sum_{j=0}^{\infty}c_{j}\int_{0}^{\infty}u^{\frac{s+q}{2}+j-1}e^{-ue^{i{\tilde \theta_{k}}}}
du \cdot t^{-\frac{s+q}{2}-j} \rbrace .\qquad $$
For $q\geq 1$, the zero coefficients in the asymptotic expansions of 
$\zeta_{q}(0)$ and $- \zeta_{q}^{\prime}(0)$ as $t\rightarrow\infty$
can be obtained only from the term
$$\frac{\Gamma(s+q)}{\Gamma(\frac{s+q}{2})\Gamma(s)} a_{m-1}
\int_{0}^{\infty}u^{\frac{s}{2}-1}e^{-ue^{i{\tilde \theta_{k}}}} du \cdot t^{-\frac{s}{2}}. \tag3.6 $$
One can check by using (3.3) that
$$
\int_{0}^{\infty} (ue^{i{\tilde \theta_{k}}})^{\frac{s}{2}-1}
e^{-ue^{i{\tilde \theta_{k}}}}(e^{i{\tilde \theta_{k}}}) du = \int_{0}^{\infty}r^{\frac{s}{2}-1}e^{-r}dr
=\Gamma (\frac{s}{2}).$$
Hence, the equation (3.6) can be simplied as
$$
\frac{\Gamma(s+q)\Gamma (\frac{s}{2})}{\Gamma(\frac{s+q}{2})\Gamma(s)} 
e^{-i\frac{{\tilde \theta_{k}}}{2}s} a_{m-1} t^{-\frac{s}{2}}.
\tag3.7 $$

\proclaim{Lemma 3.1}
For each $q\geq 1$, $a_{m-1}=0$. Therefore, the zero coefficients in the asymptotic expansions of
$\zeta_{q}(0)$ and $- \zeta_{q}^{\prime}(0)$ as $t\rightarrow\infty$ are zero.
\endproclaim
{\it Proof} \hskip 0.3 true cm
We denote by $\beta(s)$
$$\beta(s)=\frac{1}{\Gamma(s)}\int_{0}^{\infty}r^{\frac{s+q}{2}-1}
Tr\lbrace \sqrt{e^{-i\phi_{k}}B^{2}}^{q}e^{-re^{-i\phi_{k}}B^{2}} \rbrace dr . \qquad \qquad \qquad $$
Then from (3.5),
$$\beta(s)=\frac{1}{\Gamma(s)}\int_{0}^{1}r^{\frac{s+q}{2}-1}Tr\lbrace \sqrt{e^{-i\phi_{k}}B^{2}}^{q}
e^{-re^{-i\phi_{k}}B^{2}} \rbrace dr \qquad \qquad \qquad \qquad $$
$$ \qquad \qquad  + \quad \frac{1}{\Gamma(s)}\int_{1}^{\infty}r^{\frac{s+q}{2}-1}
Tr\lbrace \sqrt{e^{-i\phi_{k}}B^{2}}^{q}e^{-re^{-i\phi_{k}}B^{2}} \rbrace dr $$
$$=\sum_{j=0}^{N}a_{j} \frac{1}{\Gamma(s)}\int_{0}^{1}r^{\frac{s+q}{2}-1} r^{\frac{j-(m-1)-q}{2}} dr +
\frac{1}{\Gamma(s)}\Psi_{N}(s) \qquad \qquad \qquad  $$ 
$$ \quad = \sum_{j=0}^{N}a_{j}\frac{1}{\Gamma(s)}\frac{2}{s+j-(m-1)} + \frac{1}{\Gamma(s)} \Psi_{N}(s), 
\qquad \qquad \qquad \qquad \qquad \quad $$
where $\Psi_{N}(s)$ is holomorphic for $Res > -q$.
Hence, $$\beta(0)=2a_{m-1}. \tag3.8 $$
On the other hand,
$$\beta(s)=\sum_{\lambda_{j}\in Spec(|B|)} \frac{1}{\Gamma(s)}\int_{0}^{\infty}r^{\frac{s+q}{2}-1}
\sqrt{e^{-i\phi_{k}}\lambda_{j}^{2}}^{q}e^{-re^{-i\phi_{k}}\lambda_{j}^{2}} dr.\qquad $$
One can check by contour integration ({\it c.f.} (3.3)) that
$$
\beta(s)=\frac{\Gamma(\frac{s+q}{2})}{\Gamma(s)}\zeta_{e^{-i\frac{\phi_{k}}{2}}|B|}(s)
=e^{i\frac{\phi_{k}}{2}s}\frac{\Gamma(\frac{s+q}{2})}{\Gamma(s)}\zeta_{|B|}(s). \qquad \qquad \qquad $$
Since $\zeta_{|B|}(s)$ is regular at $s=0$, $\beta(0)=0$ and
hence from (3.8) $a_{m-1}=0$.
\qed

\vskip 0.3 true cm

From Lemma 3.1, it's enough to consider $\zeta_{0}(s)$ to compute the zero coefficients in the asymptotic expansions of 
$\zeta(0)$ and $- \zeta^{\prime}(0)$ as $t\rightarrow\infty$. Recall that from (3.4)
$$
\zeta_{0}(s)= \frac{1}{\Gamma(\frac{s}{2})} \hskip 0.1 true cm  t^{-\frac{s}{2}}
\int_{0}^{\infty} u^{\frac{s}{2}-1} e^{-ue^{i{\tilde \theta_{k}}}}
Tr\lbrace e^{-\frac{u}{t}e^{-i\phi_{k}}B^{2}} \rbrace du $$
Then, as $t\rightarrow\infty$,
$$\zeta_{0}(s) \sim \sum_{j=0}^{\infty} \hskip 0.1 true cm d_{j} \frac{1}{\Gamma(\frac{s}{2})} \hskip 0.1 true cm  t^{-\frac{s}{2}}
\int_{0}^{\infty} u^{\frac{s}{2}-1} e^{-ue^{i{\tilde \theta_{k}}}} 
\left(\frac{u}{t}\right)^{\frac{j-(m-1)}{2}}du \qquad $$
$$ = \sum_{j=0}^{\infty} \hskip 0.1 true cm  d_{j} \frac{1}{\Gamma(\frac{s}{2})} \hskip 0.1 true cm  t^{-\frac{s+j-(m-1)}{2}}
\int_{0}^{\infty} u^{\frac{s+j-(m-1)}{2}-1} e^{-ue^{i{\tilde \theta_{k}}}} du $$
$$= \sum_{j=0}^{\infty} \hskip 0.1 true cm  d_{j} \left( e^{-i{\tilde \theta_{k}}}\right)^{\frac{s+j-(m-1)}{2}}
\frac{\Gamma(\frac{s+j-(m-1)}{2})}{\Gamma(\frac{s}{2})} \hskip 0.1 true cm  t^{-\frac{s+j-(m-1)}{2}}. \tag3.9 $$
Hence, the zero coefficients $\pi_{0}(\zeta_{0}(0))$, $\pi_{0}(-\zeta_{0}^{\prime}(0))$ in the asymptotic expansions of 
$\zeta_{0}(0)$ and$ - \zeta_{0}^{\prime}(0)$ as $t\rightarrow\infty$ come from only the term 
$$d_{m-1}\left( e^{-i{\tilde \theta_{k}}}\right)^{\frac{s}{2}}t^{-\frac{s}{2}}, \tag3.10 $$
and hence, 
$$\pi_{0}(\zeta_{0}(0))=d_{m-1}. \qquad \qquad \qquad \qquad \qquad   \tag3.11 $$
$$\pi_{0}(-\zeta_{0}^{\prime}(0))=\frac{i}{2}{\tilde \theta_{k}}d_{m-1}
= \frac{i}{2}(\theta_{k}-\phi_{k})d_{m-1}. \tag3.12$$
Therefore, from (3.1), (3.11) and (3.12)  the zero coefficient $c_{k}$ in the asymptotic expansion of 
$\log Det(\sqrt{B^{2}+\alpha_{k}t \hskip 0.1 true cm}+|B|)$ is as follows.
$$c_{k}=\frac{i\phi_{k}}{2}d_{m-1}+\frac{i}{2}(\theta_{k}-\phi_{k})d_{m-1}
= \frac{i}{2}\theta_{k} d_{m-1}.$$
Since 
$$d_{m-1}=\zeta_{e^{-i\phi_{k}}B^{2}}(0)+dimKer ( e^{-i\phi_{k}}B^{2})
= \zeta_{B^{2}}(0)+dimKer B^{2},  $$
$d_{m-1}$ does not depend on $e^{-i\phi_{k}}$ and hence
$$c=\frac{1}{m}\sum_{k}c_{k}=\frac{i}{2m}d_{m-1}\sum_{k}\theta_{k}=0,\qquad $$
which completes the proof of Theorem 1.4.

\vskip 1 true cm

\S 4 {\bf Proof of Theorem 1.5 }

\vskip 0.5 true cm

In this section we are going to prove Theorem 1.5. We first assume that both $Q_{i}+\sqrt{B^{2}}$ ($i=1$, $2$) are invertible operators. 
Then this implies that $KerB=0$  \newline ( Proposition 4.5).
From Theorem 1.2, 1.3 and 1.4 we have :
$$ 
\log Det {\frak D}_{M_{r}}^{2} -  \log Det {\frak D}_{M_{1,r},P_{<}}^{2} - \log Det {\frak D}_{M_{2,r},P_{>}}^{2}  $$
\vskip -2mm
$$ =
\log Det {\frak D}_{M_{r}}^{2} - \log Det {\frak D}_{M_{1,r},D_{0}}^{2} - \log Det {\frak D}_{M_{2,r},D_{0}}^{2} + \qquad \qquad \qquad \qquad \qquad \qquad $$
$$ \qquad \log Det {\frak D}_{M_{1,r},D_{0}}^{2} + \log Det {\frak D}_{M_{2,r},D_{0}}^{2} 
-  \log Det {\frak D}_{M_{1,r},P_{<}}^{2} - \log Det {\frak D}_{M_{2,r},P_{>}}^{2}  $$
\vskip -2mm
$$
= \log Det {\frak D}_{M_{r}}^{2} - \log Det {\frak D}_{M_{1,r},D_{0}}^{2} - \log Det {\frak D}_{M_{2,r},D_{0}}^{2}  \qquad \qquad \qquad \qquad \qquad \qquad \quad $$
$$ 
+ \log Det (-\partial_{u}^{2}+B^{2})_{N_{-r,0},D_{-r},D_{0}} + \log Det (-\partial_{u}^{2}+B^{2})_{N_{0,r},D_{0},D_{r}} \quad $$
$$ - \log Det (-\partial_{u}^{2}+B^{2})_{N_{-r,0},D_{-r},P_{<}} - \log Det (-\partial_{u}^{2}+B^{2})_{N_{0,r},P_{>},D_{r}}  \quad $$
$$ \qquad + \log Det R_{M_{1,r},D_{-r}} + \log Det R_{M_{2,r},D_{r}} - 
\log Det R_{M_{1,r},P_{<}} - \log Det R_{M_{2,r},P_{>}}
$$
\vskip -2mm
$$= \log Det {\frak D}_{M_{r}}^{2} - \log Det {\frak D}_{M_{1,r},D_{0}}^{2} - \log Det {\frak D}_{M_{2,r},D_{0}}^{2}  
 \qquad \qquad \qquad \qquad \qquad \qquad \quad $$
$$ \qquad + \log Det R_{M_{1,r},D_{-r}} + \log Det R_{M_{2,r},D_{r}} - 
\log Det R_{M_{1,r},P_{<}} - \log Det R_{M_{2,r},P_{>}} $$
$$ - \log Det Q_{(\partial_{u}+|B|),r} . \qquad \qquad \qquad \qquad \qquad \qquad \quad
\qquad \qquad \quad $$
From Theorem 1.1, we have :
$$ \lim_{r \to \infty} \lbrace \log Det {\frak D}_{M_{r}}^{2} -  \log Det {\frak D}_{M_{1,r},P_{<}}^{2} - 
\log Det {\frak D}_{M_{2,r},P_{>}}^{2} \rbrace  $$
\vskip -2mm
$$= \frac{1}{2}\log Det (B^{2}) + 
\lim_{r \to \infty} \lbrace - \log Det Q_{(\partial_{u}+|B|),r}  + \log Det R_{M_{1,r},D_{-r}}
\quad \qquad \qquad \qquad \qquad \qquad$$ 
$$ + \log Det R_{M_{2,r},D_{r}} - \log Det R_{M_{1,r},P_{<}} - \log Det R_{M_{2,r},P_{>}} \rbrace \tag4.1 $$

\vskip 0.3 true cm

Now we describe the operators $Q_{(\partial_{u}+|B|),r}$, $R_{M_{1,r},D_{-r}}$,  $R_{M_{2,r},D_{r}}$, $R_{M_{1,r},P_{<}}$ and  $R_{M_{2,r},P_{>}}$
in terms of $Q_{i}$ and $B$. One can check the following lemma by direct computation.
\proclaim{Lemma 4.1}
Suppose that $f\in C^{\infty}(Y)$ with $Bf= \lambda f$ and $Ker B=0$. Then :
$$(1) \quad R_{M_{1,r},D_{-r}}(f)= Q_{1}(f)+
\left( |\lambda|+\frac{2|\lambda|e^{-|\lambda|r}}{e^{|\lambda|r}-e^{-|\lambda|r}} \right) f .
\qquad \qquad \qquad \qquad \qquad \qquad $$
$$(2) \quad R_{M_{2,r},D_{r}}(f)= Q_{2} (f)+
\left( |\lambda|+\frac{2|\lambda|e^{-|\lambda|r}}{e^{|\lambda|r}-e^{-|\lambda|r}} \right) f .
\qquad \qquad \qquad \qquad \qquad \qquad $$
$$(3) \quad R_{M_{1,r},P_{<}}(f)= \cases Q_{1}(f)+
\left( |\lambda|+\frac{2|\lambda|e^{-|\lambda|r}}{e^{|\lambda|r}-e^{-|\lambda|r}} \right) f & \text{ for $\lambda < 0 $ } \\
Q_{1}(f)+|\lambda|f & \text{ for $\lambda > 0 $ } . \endcases \qquad \qquad \qquad $$
$$(4) \quad R_{M_{2,r},P_{>}}(f)= \cases Q_{2}(f)+|\lambda|f & \text{ for $\lambda < 0 $ } \\
Q_{2}(f)+\left( |\lambda|+\frac{2|\lambda|e^{-|\lambda|r}}{e^{|\lambda|r}-e^{-|\lambda|r}} \right) f & \text{ for $\lambda > 0 $ }.
\endcases \qquad \qquad \qquad $$
$$(5) \quad Q_{(\partial_{u}+|B|),r}(f)=
\left( 2|\lambda|+\frac{2|\lambda|e^{-|\lambda|r}}{e^{|\lambda|r}-e^{-|\lambda|r}} \right) f .
\qquad \qquad \qquad \qquad \qquad \qquad \qquad \quad $$
\endproclaim
\noindent
The following lemma can be also checked easily.
\proclaim{Lemma 4.2}
Let $A$ be an invertible elliptic operator of order $>0$ and $K_{r}$ be a one-parameter family of 
trace class operators
such that $\lim_{r\to\infty}Tr(K_{r})=0$.
Then
$$\lim_{r\to\infty}\log Det(A+K_{r})=\log DetA.$$
\endproclaim

\vskip 0.3 true cm

\noindent
Applying Lemma 4.2, we have :
$$ \lim_{r \to \infty} \log Det R_{M_{1,r},D_{-r}} = \lim_{r \to \infty} \log Det R_{M_{1,r},P_{<}} = 
\log Det(Q_{1}+|B|). \tag4.2 $$
$$ \lim_{r \to \infty} \log Det R_{M_{2,r},D_{r}} = \lim_{r \to \infty} \log Det R_{M_{2,r},P_{>}} = \log Det(Q_{2}+|B|). \tag4.3 $$
Note that
$$\zeta_{Q_{(\partial_{u}+|B|),r}}(s)= \sum_{\lambda\in Spec(B)}
\left( 2|\lambda|+\frac{2|\lambda|e^{-|\lambda|r}}{e^{|\lambda|r}-e^{-|\lambda|r}} \right)^{-s}.$$
From Lemma 4.2 again, we have :
$$\lim_{r \to \infty} \lbrace  \log Det Q_{(\partial_{u}+|B|),r}  \rbrace =
\log 2 \cdot  \zeta_{B^{2}}(0) + \frac{1}{2}\log Det B^{2}. \tag4.4 $$
Therefore from (4.1), (4.2), (4.3) and (4.4), we have
$$\lim_{r \to \infty} \lbrace \log Det {\frak D}_{M_{r}}^{2} -  \log Det {\frak D}_{M_{1,r},P_{<}}^{2} - 
\log Det {\frak D}_{M_{2,r},P_{>}}^{2} \rbrace = -\log 2  \cdot \zeta_{B^{2}}(0).$$
This completes the proof of Theorem 1.5.

\vskip 0.3 true cm

Finally, we are going to discuss the invertibility conditions of both $Q_{1}+\sqrt{B^{2}}$ and $Q_{2}+\sqrt{B^{2}}$. 
For this purpose we need Green's formula of the following form ({\it c.f.} Lemma 3.1 in [5]).
\proclaim{Lemma 4.3}
Let $\phi$ and $\psi$ be smooth sections on $M_{j}$ ($j=1$, $2$). Then,
$$\langle {\frak D}_{M_{j}} \phi, \psi \rangle_{M_{j}} -  \langle  \phi, {\frak D}_{M_{j}} \psi \rangle_{M_{j}}
= \epsilon_{j} \langle \phi|_{Y}, G(\psi|_{Y}) \rangle_{Y}, $$
where $\epsilon_{j}=1$ for $j=2$ and  $\epsilon_{j}=-1$ for $j=1$.
\endproclaim
\noindent
Suppose that  for $f\in C^{\infty}(Y)$, $\phi_{j}$ is the solution of ${\frak D}^{2}_{M_{j}}$ with $\phi_{j}|_{Y}=f$.
Then by Lemma 4.3
$$
\langle (Q_{1}+|B|)f, f\rangle_{Y} = \langle {\frak D}_{M_{1}}\phi_{1}, {\frak D}_{M_{1}}\phi_{1}\rangle_{M_{1}} + 
\langle (|B|-B)f,f\rangle_{Y}, \tag4.5 $$
$$
\langle (Q_{2}+|B|)f, f\rangle_{Y} = \langle {\frak D}_{M_{2}}\phi_{2}, {\frak D}_{M_{2}}\phi_{2}\rangle_{M_{2}} + 
\langle (|B|+B)f,f\rangle_{Y}. \tag4.6 $$
Hence we have  
$$f\in Ker(Q_{1}+|B|) \quad \text{ if and only if } \quad {\frak D}_{M_{1}}\phi_{1}=0 \text{ and } f\in Im P_{\geq}. \tag4.7 $$
$$f\in Ker(Q_{2}+|B|) \quad \text{ if and only if } \quad {\frak D}_{M_{2}}\phi_{2}=0 \text{ and } f\in Im P_{\leq}. \tag4.8 $$
$\phi_{1}$ satisfying (4.7) can be expressed on the cylinder part by
$$\phi_{1}=\sum_{j=1}^{k}a_{j}g_{j}+\sum_{\lambda_{j}>0}b_{j}e^{-\lambda_{j}u}h_{j},  \tag4.9 $$
where $Bg_{j}=0$, $Bh_{j}=\lambda_{j}h_{j}$ and $k=dimKer B$. $\phi_{2}$ satisfying (4.8) can be expressed in the similar way.

\proclaim{Definition 4.4}
We denote $M_{1,\infty}:=M_{1} \cup_{\partial M_{1}} Y \times [0,\infty)$ and by ${\frak D}_{M_{1,\infty}}$, $E_{1,\infty}$
the natural extensions of ${\frak D}_{M_{1,r}}$, $E_{1,r}$ to $M_{1,\infty}$.
A section $\psi$ of  $E_{1,\infty}$ is called an extended $L^{2}$-solution of ${\frak D}_{M_{1,\infty}}$ if $\psi$
is a solution of ${\frak D}_{M_{1,\infty}}$ which takes
the form (4.9) on the cylinder part $[0, \infty)\times Y$. In this case $\sum_{j=1}^{k}a_{j}g_{j}$ is called 
the limiting value of $\psi$. We can define the same notions for ${\frak D}_{M_{2,\infty}}$ on $M_{2,\infty}$.
\endproclaim

It is a well-known fact that if $L$ is the set of all limiting values of extended 
$L^{2}$-solutions of ${\frak D}_{M_{1,\infty}}$, $L$ is a Lagrangian subspace
of $Ker B$ and in particular $dim L = \frac{1}{2} dim Ker B$ ({\it c.f.} [1], [2], [5], [10]).
From Definition 4.4,
$\phi_{1}$ satisfying (4.7) is the restriction of an extended $L^{2}$-solution of ${\frak D}_{M_{1,\infty}}$ on 
$M_{1,\infty}$.
We can say similar assertion for $\phi_{2}$ and
as a consequence we have the following proposition.
\proclaim{Proposition 4.5}
The invertibility of $Q_{1}+\sqrt{B^{2}}$ and $Q_{2}+\sqrt{B^{2}}$ is equivalent to the non-existence of extended $L^{2}$-solutions
of  ${\frak D}_{M_{1,\infty}}$ and ${\frak D}_{M_{2,\infty}}$ on $M_{1,\infty}$ and $M_{2,\infty}$.
Furthermore, this condition implies the invertibility of $B$,  ${\frak D}^{2}_{M_{1,r},P_{<}}$,  ${\frak D}^{2}_{M_{2,r},P_{>}}$ for each $r>0$ and
the invertibility of ${\frak D}^{2}_{M_{r}}$ for $r$ large enough.
\endproclaim
{\it Proof} \hskip 0.3 true cm
We need to prove the second assertion. $Ker B=0$ implies the invertibility of $B$.
Suppose that $\psi\in Ker{\frak D}^{2}_{M_{1,r},P_{<}}$. 
By Lemma 4.3 and (1.2) we have
$$0= \langle {\frak D}^{2}_{M_{1,r}}\psi, \psi \rangle_{M_{1,r}}=
\langle {\frak D}_{M_{1,r}}\psi, {\frak D}_{M_{1,r}}\psi \rangle_{M_{1,r}} -
\langle (\partial_{u}+B)\psi|_{Y_{0}}, \psi|_{Y_{0}} \rangle_{Y_{0}} $$
$$= \langle {\frak D}_{M_{1,r}}\psi, {\frak D}_{M_{1,r}}\psi \rangle_{M_{1,r}} .$$
Hence, ${\frak D}_{M_{1,r}}\psi=0$ and by (4.7) $\psi|_{Y_{0}}$ belongs to $Ker (Q_{1}+|B|)$, which implies that
$\psi=0$ and hence $Ker{\frak D}_{M_{1,r},P_{<}}=0$.
Since ${\frak D}_{M_{1,r},P_{<}}$ is a self-adjoint Fredholm operator ({\it c.f.} Proposition 2.4 in [5]), ${\frak D}_{M_{1,r},P_{<}}$ is an
invertible operator and so is ${\frak D}^{2}_{M_{1,r},P_{<}}$. 
We can use the same argument for ${\frak D}^{2}_{M_{2,r},P_{>}}$.

Now let us consider ${\frak D}^{2}_{M_{r}}$. 
Since ${\frak D}_{M_{r}}$ is a self-adjoint Fredholm operator, it's enough to show that $Ker {\frak D}_{M_{r}} = 0$.
From the decomposition 
$M_{r}=(M_{1}\cup M_{2})\cup_{(\partial M_{1} \cup \partial M_{2})} N_{-r,r}$, we define the Dirichlet-to-Neumann operator
$$ R_{-r,r} : C^{\infty}(Y_{-r}) \oplus C^{\infty}(Y_{r}) \rightarrow C^{\infty}(Y_{-r}) \oplus C^{\infty}(Y_{r}) $$
as follows. 
For $(f,g)\in C^{\infty}(Y_{-r}) \oplus C^{\infty}(Y_{r})$, choose $\phi_{i}\in C^{\infty}(M_{i})$ ($i=1, 2$), $\psi\in C^{\infty}(N_{-r,r})$
such that
$$
{\frak D}^{2}_{M_{i}}\phi_{i}=0, \qquad  (-\partial_{u}^{2}+B^{2})\psi=0, $$
$$ \phi_{1}|_{\partial M_{1}}= \psi|_{Y_{-r}} = f, \qquad  \phi_{2}|_{\partial M_{2}}= \psi|_{Y_{r}} =g .$$
Then we define
$$ \split
R_{-r,r}(f,g) & = \left( (\partial_{u}\phi_{1})|_{Y_{-r}} - (\partial_{u}\psi)|_{Y_{-r}}, 
- (\partial_{u}\phi_{2})|_{Y_{r}} + (\partial_{u}\psi)|_{Y_{r}} \right)    \\
& = \left( Q_{1}f - (\partial_{u}\psi)|_{Y_{-r}},  Q_{2}g + (\partial_{u}\psi)|_{Y_{r}} \right) .
\endsplit $$
If $(f,g)\in Ker R_{-r,r}$, then $\phi_{1} \cup_{Y_{-r}} \psi \cup_{Y_{r}} \phi_{2}$ is 
a smooth section which belongs to $Ker{\frak D}^{2}_{M_{r}}$ and vice versa.
Hence $Ker R_{-r,r}=0$ if and only if $Ker {\frak D}^{2}_{M_{r}} =Ker {\frak D}_{M_{r}} =0$.
By direct computation ({\it c.f.} [9]), one can check that

$$R_{-r,r} \binom{f}{g} = \left( \matrix
Q_{1}+ |B| & 0 \\ 0 & Q_{2}+ |B|
\endmatrix \right) \binom{f}{g}
+ A_{r}
\left( \matrix e^{-2r|B|} & - 1 \\ - 1 & e^{-2r |B|}  \endmatrix \right) \binom{f}{g}, $$
\noindent
where $A_{r}=\frac{2|B|}{e^{2r|B|}-e^{-2r|B|}}$. Then we have
$$
\left< R_{-r,r}\binom{f}{g}, \binom{f}{g}\right>_{L^{2}(Y)} $$
$$ \multline 
=\langle (Q_{1}+|B|)f,f\rangle + \langle (Q_{2}+|B|)g,g\rangle  \\
+ \langle A_{r}e^{-2r|B|}f,f\rangle + \langle A_{r}e^{-2r |B|}g,g\rangle
- \langle A_{r}g,f\rangle- \langle A_{r}f,g\rangle. \endmultline$$
Note that each $Q_{i}+|B|$ is a non-negative operator by (4.5), (4.6).
Let $\lambda_{0}$ be the minimum of the eigenvalues of $Q_{1}+\sqrt{\Delta_{Y}}$ and $Q_{2}+\sqrt{\Delta_{Y}}$.
Since $\lim_{r\to\infty}||A_{r}||_{L^{2}}=0$, one can choose $r_{0}$ so that for $r\geq r_{0}$,  $||A_{r}||_{L^{2}} < \lambda_{0}$. 
Then for $r\geq r_{0}$, $R_{-r,r}$ is injective and this completes the proof of Proposition 4.5.
\qed

\vskip 1 true cm

\vskip 0.3 true cm
E-mail address : ywlee\@math.inha.ac.kr

\enddocument